\documentclass[twocolumn]{autart_edited}    

\usepackage{graphicx}          

\usepackage{cite}
\usepackage{amsmath}
\usepackage{amssymb}
\usepackage{algpseudocode}
\usepackage{bbm}
\usepackage{color}
\usepackage{subfigure}

\newtheorem{definition}{Definition}
\newtheorem{proposition}{Proposition}
\newtheorem{remark}{Remark}

\newtheorem{theorem}{Theorem}

\newcommand{\Int}{\mathsf{int}}

\newcommand{\D}{\mathcal{D}}
\newcommand{\M}{\mathcal{M}}
\newcommand{\A}{\mathcal{A}}
\newcommand{\R}{\mathcal{R}}
\newcommand{\DMM}{\left[\partial \mathcal{M}\right]_-}
\newcommand{\DMO}{\left[\partial \mathcal{M}\right]_0}
\newcommand{\Mcomp}{\mathcal{M}^{\textsf{C}}}
\newcommand{\RR}{\mathbb{R}}

\begin{document}

\begin{frontmatter}
\runtitle{On Maximal Robust Positively Invariant Sets in Constrained Nonlinear Systems}  

\title{On Maximal Robust Positively Invariant Sets in Constrained Nonlinear Systems\thanksref{footnoteinfo}} 

\thanks[footnoteinfo]{© 2020. This manuscript version is made available under the CC-BY-NC-ND 4.0 license http://creativecommons.org/licenses/by-nc-nd/4.0/. DOI: https://doi.org/10.1016/j.automatica.2020.109044. Corresponding author: Willem Esterhuizen.}

\author[ACSD]{Willem Esterhuizen}\ead{willem.esterhuizen@etit.tu-chemnitz.de},
\author[ACSD]{Tim Aschenbruck}\ead{tim.aschenbruck@etit.tu-chemnitz.de},
\author[ACSD]{Stefan Streif}\ead{stefan.streif@etit.tu-chemnitz.de}

\address[ACSD]{Technische Universit\"{a}t Chemnitz, Automatic Control and System Dynamics Laboratory, Germany}     
          
\begin{keyword}                           
	invariant sets, constraint satisfaction problems, nonlinear systems
\end{keyword}                            

\begin{abstract}       
	In this technical communiqu\'{e} we study the maximal robust positively invariant set for state-constrained continuous-time nonlinear systems subjected to a
	bounded disturbance. Extending results from the theory of barriers, we show that this set is closed and that its boundary consists of two complementary parts, one of which we name the invariance barrier, which consists of trajectories that satisfy the maximum principle.
\end{abstract}

\end{frontmatter}

\section{Introduction}
Set invariance is a fundamental concept in control theory due to its well-known relationship with stability, \cite{BLANCHINI_1999,Blanchini:2015}. This paper focusses on the \emph{maximal robust positively invariant set} (MRPI) of a continuous-time nonlinear system subjected to state constraints and a bounded disturbance term: any trajectory initiating in this set remains in it for all future time, regardless of the disturbance realisation. Robust invariant sets play a central role in the design of obstacle-avoiding path-planning methodologies, \cite{Blanchini_2004,Danielson_2016}; in the design of reference governors, \cite{Kolmanovsky_2014}; and in stability and recursive feasibility studies in predictive control, \cite{Kerrigan_2000,Mayne2000}. The MRPI has also recently been applied to the study of epidemics, \cite{ESTERHUIZEN_2020c}. Most work on the MRPI focusses on its approximation via algorithms that involve the iterative computation and intersection of one-step \emph{predecessor sets}\footnote{In discrete-time, this is the set of all states for which the subsequent state is contained in a set $S$, for all disturbance inputs. These sets go under various names in the literature.} \cite{Rakovic_2008,Cannon_2015,Kolmanovsky_1998}. Other sets that are closely related to the MRPI, but should not be confused with it, include: \textit{minimal robust positively invariant sets}, \cite{Kerrigan_2005}; \emph{regions of attraction}, \cite{Henrion_2014_b}; and \emph{backwards reachable sets}, \cite{Mitchel2005}.

In this paper we characterise the MRPI's boundary: we adapt results from the \emph{theory of barriers}, \cite{DeDona_siam,ESTERHUIZEN_2016}, and show that many facts concerning the so-called \emph{admissible set} carry over to the MRPI analogously. Under some assumptions, we show that the MRPI is closed and that its boundary consist of two complementary parts: 
one called the \emph{usable part}, which is contained in the constraint set's boundary, the other, which we name the \emph{invariance barrier}, contained in the constraint set's interior. Furthermore, we show that the invariance barrier may consist of parts that are made up of special integral curves of the system that satisfy Pontryagin's maximum principle, a fact that may be used to construct MRPIs. We emphasise that these curves give an \emph{exact} description of the invariance barrier, and if analytic solutions to them are not obtainable, the error in the computation of the set would depend on the integration scheme used.

\vspace{-0.21cm}
The outline of the paper is as follows. In Section~\ref{sec_problem_formulation} we present the constrained system under study. In Section~\ref{sec_MRPI_closedness} we show that the MRPI is closed, and in Section~\ref{sec_boundary} we characterise its boundary. Section~\ref{sec_ult_tang} is dedicated to the ultimate tangentiality condition, which is satisfied at the intersection of the invariance barrier and the boundary of the constrained state-space. Section~\ref{sec_main_result} presents our main result: Theorem~1. Section~\ref{sec_examples} contains examples and Section~\ref{sec_conclusion} concludes the paper.

\section{Constrained System Formulation}\label{sec_problem_formulation}

We consider the following nonlinear system:
\begin{align}
\dot{x}(t) & = f(x(t),d(t)),\,\,x(t_0) = x_0,\,\,d \in \mathcal{D},\label{sys_eq_1}\\
g_i(x(t)) &\leq 0,\forall t\in[t_0,\infty),\,\,i=1,2,\dots,p,\label{sys_eq_4}
\end{align}
where $x(t)\in\mathbb{R}^n$ is the state and $d(t)\in\mathbb{R}^m$ is a disturbance input. We make similar assumptions to \cite{DeDona_siam}:
\begin{description}
	\item[(A1)] The space $\mathcal{D}$ is the set of all Lebesgue measurable functions that map the interval $[t_0,\infty)$ to a set $D\subset\RR^m$, which is compact and convex.
	\item[(A2)] The function $f$ is $C^2$ with respect to $d\in D$, and for every $d$ in an open subset containing $D$, the function $f$ is $C^2$ with respect to $x\in\RR^n$.
	\item[(A3)] Every $x_0$, with $\Vert x_0 \Vert < \infty$, and every $d\in \mathcal{D}$ admits a unique absolutely continuous integral curve of \eqref{sys_eq_1} that remains bounded over any finite time interval.
	\item[(A4)] The set $f(x,D)\triangleq\{f(x,d) : d\in D\}$ is convex for all $x\in \RR^{n}$.
	\item[(A5)] For every $i$ the function $g_i$ is $C^{2}$ with respect to $x\in\RR^n$, and the set $\{x:g_i(x) = 0\}$ defines a manifold.
\end{description}
The assumptions (A2)-(A5) are needed to arrive at a compactness result, stated in Proposition~\ref{prop_1_compactness}, which is taken from \cite{DeDona_siam}. By $x^{(d,x_0,t_0)}$ we will refer to the solution of \eqref{sys_eq_1} with the initial condition $x_0\in\RR^n$ at time $t_0\in\mathbb{R}$ and a disturbance realisation $d\in\mathcal{D}$. If the initial time is clear from context we will use the notation $x^{(d,x_0)}$, and if the initial condition is clear we will use $x^{(d)}$. By $x^{(d,x_0,t_0)}(t)$, $x^{(d,x_0)}(t)$ and $x^{(d)}(t)$, with $t\in[t_0,\infty)$, we will refer to the solution at time $t$. Given two disturbance realisations, $d_1\in\mathcal{D}$ and $d_2\in\mathcal{D}$, along with a time instant $\tau\in[t_0,\infty)$, the concatenated disturbance given by $d_3(t) = 
\begin{cases} 
d_1(t)\,\,\text{for}\,\,t\in[t_0,\tau)\\
d_2(t)\,\,\text{for}\,\,t\in[\tau,\infty)
\end{cases}$ also satisfies $d_3\in\mathcal{D}$. We denote this concatenation by $d_3 = d_1\bowtie_{\tau} d_2$. Let $g(x)\triangleq(g_1(x),g_2(x),\dots,g_p(x))^T$. 
By $\mathbb{I}(x)$ we refer to the set $\{i\in\{1,2,\dots,p\} : g_i(x) = 0\}$. We introduce the following sets in order to lighten our notation: $G\triangleq\{x : g_i(x)\leq 0, i=1,\dots,p\}$, $G_-\triangleq\{x : g_i(x) < 0, i=1,\dots,p\}$, $G_0\triangleq\{x : \exists i, g_i(x) = 0\}$. The notation $L_fg(x, d)\triangleq \nabla g(x).f(x,d)$ denotes the Lie derivative of a differentiable function $g$ with respect to the vector field $f(x,d)$ at the point $x$. If $S$ is a set, then $\Int(S)$ denotes its interior, $\mathsf{cl}(S)$ its closure and $S^\mathsf{C}$ its complement.

\section{Some properties of the MRPI}\label{sec_MRPI_closedness}
We recall some notions from the literature on invariant sets, see for example \cite{BLANCHINI_1999,Blanchini:2015,Kerrigan_2000, Henrion_2019}.
\begin{definition}
	A set $\Omega\subset\mathbb{R}^n$ is a \emph{robust positively invariant set} (RPI) of the system \eqref{sys_eq_1} provided that $x^{(d,x_0,t_0)}(t)\in\Omega$ for all $t\in[t_0,\infty)$, for all $x_0\in\Omega$ and for all $d\in\mathcal{D}$.
\end{definition}
\begin{definition}\label{def:MRPI}
	We denote by $\M$ the \emph{maximal robust positively invariant set} (MRPI) of the system \eqref{sys_eq_1}-\eqref{sys_eq_4} contained in $G$. In other words, $\M$ is the union of all RPIs that are subsets of $G$.
\end{definition}
Next, we introduce an equivalent expression of the MRPI which will make it easier to study. 
\begin{proposition}\label{prop_1}
	An equivalent definition of $\mathcal{M}$ for system \eqref{sys_eq_1}-\eqref{sys_eq_4} is given by:
	\begin{equation}
	\R =\{x_0\in G : x^{(d,x_0,t_0)}(t)\in G,\,\,\forall t\in[t_0,\infty),\,\,\forall d\in\mathcal{D}\}.\nonumber
	\end{equation}
	In other words, $\mathcal{M} = \R$.
\end{proposition}
\vspace{-0.7cm}
\begin{pf}
	If $\R$ were not an RPI there would exist an $x_1\in\R$, $d_1\in\mathcal{D}$, and $t_1\in[t_0,\infty)$ such that $x^{(d_1,x_1,t_0)}(t_1)\triangleq x_2 \notin \R$. Then, there would exist a $d_2\in\D$ and $t_2\in[t_1,\infty)$ such that $x^{(d_2,x_2,t_1)}(t_2)\notin G$. We could then form the disturbance $d_3 = d_1\bowtie_{t_1}d_2$ for which $x^{(d_3,x_1,t_0)}(t_2)\notin G$, contradicting the fact that $x_1\in\R$. Clearly $\R\subset G$, and so $\R\subset\M$. For any $\hat{x}\in\M$, $x^{(d,\hat{x},t_0)}(t)\in\mathcal{M}\subset G$ for all $t\in[t_0,\infty)$, for all $d\in\mathcal{D}$, by Definition~\ref{def:MRPI}, thus $\mathcal{M}\subset\R$, and thus $\mathcal{M}=\R$.
\end{pf}
\vspace{-0.4cm}
\begin{remark}
	The \emph{admissible set} is defined as follows:
	\[
		\A\triangleq \{x_0\in G:\exists d\in\D\,\,\mathrm{s.t.}\,\,x^{(d,x_0,t_0)}(t)\in G,\forall t\in[t_0,\infty)\},
	\]
	and so the equivalent description of $\M$, as in Proposition~\ref{prop_1}, exposes the relationship between $\A$ and $\M$. Indeed, the results from \cite{DeDona_siam} adapt in an intuitive way (the ``mins'' turn into ``maxes''), but the proofs are not always obvious.
\end{remark}

\subsection{Closedness of the MRPI}

\begin{proposition}\label{prop_closed}
	Under (A1)-(A5), $\mathcal{M}$ is closed.
\end{proposition}
This results follow the same lines of reasoning as the proof of closedness of the set $\A$, see \cite[Prop. 4.1]{DeDona_siam}. Briefly, consider the compactness result, as summarised in \cite[Lem. A.2]{DeDona_siam}:
\begin{proposition}\label{prop_1_compactness}
	Assume (A1)-(A4) hold. Let $\mathcal{X}(x_0)$ be the set of all integral curves initiating from $x_0\in\RR^n$, satisfying \eqref{sys_eq_1}. Consider $\mathcal{X} = \cup_{x_0\in\mathcal{X}_0}\mathcal{X}(x_0)$, with $\mathcal{X}_0$ a compact subset of $\RR^n$. From every sequence $\{x^{(d_k,x_k)}\}_{k\in\mathbb{N}}\subset\mathcal{X}$ one can extract a uniformly convergent subsequence on
	every finite interval $[0, T]$, whose limit $\xi$ is an absolutely continuous integral curve on $[0,\infty)$, belonging to $\mathcal{X}$.
\end{proposition}
One can now consider an arbitrary disturbance realisation, $d\in\D$, along with a sequence of initial conditions, $\{x_k\}_{k\in\mathbb{N}}\subset\M$, converging to a point $\bar{x}\in\mathbb{R}$. Because there exists a subsequence $\{x^{(d,x_{k_l},t_0)}\}_{l\in\mathbb{N}}$ converging to an integral curve of the system, and because the $g_i$'s are continuous, one can argue that $\bar{x}\in\M$.

\section{The boundary of the MRPI}\label{sec_boundary}
Focussing on the set's boundary, we let $\DMO \triangleq \partial\mathcal{M}\cap G_0$ and $\DMM\triangleq\partial\mathcal{M}\cap G_-$, which we call the \emph{usable part} and \emph{invariance barrier}, respectively. Clearly, the set $\DMO$ coincides with the set of points $z\in G_0$ for which $\max_{d\in D}	\max_{i\in\mathbb{I}(z)} L_f g_i(z,d)\leq 0$. Turning our attention to $\DMM$, consider the following set, where $T<\infty$:
$
\mathcal{M}_T = \{x_0 : x^{(d,x_0,t_0)}(t)\in G,\,\,\forall t\in[t_0,T],\,\,\forall d\in\mathcal{D}\}.\label{def_R_T}
$
If we refer to $\M$ as defined in Proposition~\ref{prop_1}, then clearly $\M\subset\M_{T_1}\subset\M_{T_2}\subset G$, for $0\leq T_2\leq T_1 < \infty$. The sets $\M_{T}$ are generally not robustly invariant, but when this is the case we have $\M_T\subset\M$, and thus $\M = \M_T$.\footnote{This is related to the set being \emph{finitely-determined} in discrete-time, see \cite{Kolmanovsky_1998}.} We introduce the following assumption:
\begin{description}
	\item[(A6)] There exists a $T<\infty$ such that $\M_T = \M$.
\end{description}
\begin{proposition}\label{prop_stay_in_DMM}
	Assume (A1)-(A6) hold. Consider a point $\bar{x}\in\DMM$. Then there exists $\bar{d}\in\D$ such that the corresponding integral curve runs along $\DMM$ and intersects $G_0$ in finite time.
\end{proposition}
\vspace{-0.6cm}
\begin{pf}	
	Consider a sequence $\{x_k\}_{k\in\mathbb{N}}$, with $x_k\in\Mcomp$ for all $k$, converging to a point $\bar{x}\in\DMM$. For every $x_k$, there exists a $d_k\in\D$, a $t_k\in[t_0,T]$ and an $i_k\in\{1,2\dots,p\}$ such that $g_{i_k}(x^{(d_k,x_k,t_0)}(t_k)) > 0$. Note that each $x^{(d_k,x_k,t_0)}$ is entirely contained in $\Mcomp$. Using the compactness result from Proposition~\ref{prop_1_compactness}, we can select a subsequence $\{x^{(d_{k_l},x_{k_l},t_0)}\}_{l\in\mathbb{N}}$ that uniformly converges on the interval $[t_0,T]$ to an integral curve belonging to the system \eqref{sys_eq_1}. Moreover, we can do so such that the sequence $\{t_{k_l}\}_{{k_l}\in\mathbb{N}}$ (satisfying $t_{k_l} \leq T$) is monotonically increasing. Thus, we have constructed a sequence of integral curves that uniformly converges to the curve $x^{(\bar{d},\bar{x},t_0)}$ with $\bar{d}\in\D$, $\bar{x}\in\DMM$ and such that $g_i(x^{(\bar{d},\bar{x},t_0)}(\bar{t})) = 0$ (by the continuity of $g$) for some $\bar{t}\leq T$ and some $i\in\{1,2\dots,p\}$. From the definition of $\M$, the curve $x^{(\bar{d},\bar{x},t_0)}$ cannot intersect $\Mcomp$ for any $t\in[t_0,\bar{t}]$. Moreover, it cannot intersect the interior of $\M$ for any $t\in[t_0,\bar{t}]$, for otherwise it will not be the uniform limit of a sequence of integral curves, each entirely contained in $\Mcomp$. Thus, $x^{(\bar{d},\bar{x},t_0)}(t)\in\DMM$ for all $t\in[t_0,\bar{t})$, $x^{(\bar{d},\bar{x},t_0)}(\bar{t})\in G_0$.
\end{pf}
\begin{remark}
	If there does not exist a $T<\infty$ such that $\M_T = \M$, then an integral curve initiating on $\DMM$ may remain in $G_-$ for all future time. As an example, consider the scalar system $\dot{x}(t) = x(t) + d(t)$, with $x(t)\leq 0$ and $|d(t)|\leq 1$. Clearly $\M = (-\infty,-1]$, and for $\bar{x}\in \DMM = \{-1\}$ there does not exist a disturbance realisation such that the resulting trajectory remains in $\DMM$ and eventually intersects $G_0$.
\end{remark}

\section{Ultimate tangentiality}\label{sec_ult_tang}
The next proposition says that $\DMM$ intersects $G_0$ tangentially. We note that its proof is simpler than its analogue, \cite[Prop 6.1]{DeDona_siam}, because we can arrive at the conclusion via a contradiction argument.
\begin{proposition}\label{prop_6}
	Assume that (A1)-(A6) hold. Consider a point $\bar{x}\in\DMM$ along with a disturbance realisation, $\bar{d}\in\mathcal{D}$, as in Proposition~\ref{prop_stay_in_DMM}, such that $x^{(\bar{d},\bar{x},t_0)}(t)\in\DMM$ for all $t\in[t_0,\bar{t})$, $x^{(\bar{d},\bar{x},t_0)}(\bar{t})\triangleq z\in G_0$. Then,
	\begin{equation}
	\max_{i\in\mathbb{I}(z)} L_fg_i(z,\bar{d}(\bar{t})) = \max_{d\in D} \max_{i\in\mathbb{I}(z)}L_f g_i(z,d) = 0. \label{Prop6_eq}
	\end{equation}
\end{proposition}
\vspace{-1cm}
\begin{pf}
	\tolerance=390
	The mapping $t\rightarrow g_{i}(x^{(\bar{d},\bar{x},t_0)}(t))$ is nondecreasing for all $i\in\mathbb{I}(z)$ over an interval $(\bar{t} - \eta,\bar{t}]$ with $\eta>0$ and small enough, which implies that $\max_{i\in\mathbb{I}(z)}L_fg_i(z,\bar{d}(\bar{t}))$ $\geq 0$. Suppose that there exists $\hat{d}\in D$ such that $\max_{i\in\mathbb{I}(z)} L_fg_i(z,\hat{d}) > \max_{i\in\mathbb{I}(z)} L_fg_i(z,\bar{d}(\bar{t}))$. This would imply that $\max_{i\in\mathbb{I}(z)} L_fg_i(z,\hat{d}) > 0$, contradicting the fact that $z\in\DMO$. Thus, we have $\max_{i\in\mathbb{I}(z)} L_fg_i(z,d) \leq \max_{i\in\mathbb{I}(z)} L_fg_i(z,\bar{d}(\bar{t}))$ for all $d\in D$.
	We have established that $0\leq \max_{i\in\mathbb{I}(z)}L_fg_i(z,\bar{d}(\bar{t})) = \max_{d\in D} \max_{i\in\mathbb{I}(z)} L_fg_i(z,d) \leq 0$, to complete the proof.
\end{pf}
\vspace{-0.6cm}
\section{Construction of the invariance barrier}\label{sec_main_result}
We now present the main result. We state the theorem in full in the interest of completeness, and sketch the proof.
\begin{theorem}\label{theorem_main}
	Assume that (A1)-(A6) hold. Every integral curve $x^{(\bar{d})}$ on $\DMM$ and the corresponding disturbance realisation, $\bar{d}\in \D$, as in Proposition~\ref{prop_stay_in_DMM}, satisfy the following necessary conditions. There exists a nonzero absolutely continuous maximal solution $\lambda^{\bar{d}}$ to the adjoint equation:
	\begin{align}
	&\dot{\lambda}^{\bar{d}}(t) = -\left( \frac{\partial f}{\partial x}(x^{(\bar{d})}(t),\bar{d}(t)) \right)^T \lambda^{\bar{d}}(t),\nonumber\\
	\quad &\lambda^{\bar{d}}(\bar{t}) = (\nabla g_{i^*}(z))^T,\label{eq_adj}
	\end{align}
	\tolerance=100
	with $\bar{t}$ as time at which $x^{(\bar{d})}$ intersects $G_0$, $z\triangleq x^{(\bar{d})}(\bar{t})$, and $ L_fg_{i^*}(z,\bar{d}(\bar{t})) = \max_{i\in\mathbb{I}(z)} L_fg_i(z,\bar{d}(\bar{t}))$, such that
	\begin{equation}
	\max_{d\in D}\{\lambda^{\bar{d}}(t)^Tf(x^{(\bar{d})}(t),d) \} = \lambda^{\bar{d}}(t)^T f(x^{(\bar{d})}(t),\bar{d}(t)) = 0,\label{eq_Hamil}
	\end{equation}
	for almost every $t\in[t_0,\bar{t}]$. Moreover, at $\bar{t}$ the ultimate tangentiality condition holds:
	\begin{align}
	\max_{d\in D} \max_{i\in\mathbb{I}(z)}L_f g_i(z,d) &= \max_{i\in\mathbb{I}(z)}L_f g_i(z,\bar{d}(\bar{t}))\nonumber\\
	&= L_fg_{i^*}(z,\bar{d}(\bar{t})) = 0.\label{thm1_ult_tan}
	\end{align}
\end{theorem}
\vspace{-0.9cm}
\begin{pf}{ (Sketch)
	Consider the reachable set at time $t\in\mathbb{R}$ from a point $\bar{x}\in\mathbb{R}^n$: $X_t(\bar{x}) \triangleq \{x^{(d,\bar{x},t_0)}(t): d\in\mathcal{D}\}$. The first step in the proof's sketch is to argue that the integral curve $x^{(\bar{d})}$, that runs along $\DMM$, satisfies $x^{(\bar{d})}(t)\in\partial X_t(\bar{x})$ for all $t\in[t_0,\bar{t})$. Indeed, the proof of the analogous statement in \cite[Prop. 7.1]{DeDona_siam} adapts to the current setting, with the difference that the reachable set is now contained in $\M$. Therefore, from Theorem~\ref{thm_PMP} of the appendix, there exists a solution, labelled $\lambda^{\bar{d}}$, to the adjoint equation, \eqref{eq_1_PMP}, such that the Hamiltonian is maximised for almost every $t\in[t_0,\bar{t}]$ as in \eqref{eq_2_PMP}. We still need to show that with the final condition $\lambda^{\bar{d}}(\bar{t}) = (\nabla g_{i^*}(z))^T$ we obtain a solution to the adjoint equation such that the constant on the right-hand side of \eqref{eq_2_PMP} is zero. As described in great detail in \cite[Ch.4]{Lee_Markus} and \cite{PBGM}, $\lambda^{\bar{d}}(t)$ is the outward normal of a hyperplane that contains the elementary perturbation cone, $\mathcal{K}_t$, (see \cite[Ch.4]{Lee_Markus}) for every $t\in[t_0,\bar{t}]$. Moreover, we have:
	\begin{equation}\label{lambda_v_const}
	\lambda^{\bar{d}}(t)^Tv(t) = \lambda^{\bar{d}}(\bar{t})^Tv(\bar{t}) \leq 0\quad \forall t\in[t_0,\bar{t}],
	\end{equation}
	where $v(t)\in\mathcal{K}_t$. Let $v(t)$ be associated with the perturbation data $\{d,\tau,l\}$. Then, after dividing by $l$, we have
	$
	\lambda^{\bar{d}}(\bar{t})^T [ f(x^{(\bar{d})}(\bar{t}),d) - f(x^{(\bar{d})}(\bar{t}),\bar{d}(\bar{t})) ] \leq 0, \forall d\in D.
	$
	Recall from the proof of Proposition~\ref{prop_6} that \eqref{Prop6_eq} implies
	$
	\max_{i\in\mathbb{I}(z)} \nabla g_i(x^{(\bar{d})}(\bar{t})) [f(x^{(\bar{d})}(\bar{t}),d) - f(x^{(\bar{d})}(\bar{t}),\bar{d}(\bar{t}) ] \leq 0,$ $ \forall d\in D,$
	from where we deduce that $\lambda^{\bar{d}}(\bar{t}) = \nabla g_{i^*}(z)^T$. 
	From \eqref{lambda_v_const} and \eqref{Prop6_eq} we deduce that the constant on the right-hand side of \eqref{eq_2_PMP} is zero. The ultimate tangentiality condition was proved in Proposition~\ref{prop_6}.}
\end{pf}
\vspace{-0.5cm}
\begin{remark}
	Note that Theorem 7.1 from \cite{DeDona_siam} stated that the input function resulting in the integral curve running along the barrier of the admissible set minimises the Hamiltonian for almost every time, whereas in the case of the MRPI, it is maximised.
\end{remark}
The set $\M$ may be found with Theorem~\ref{theorem_main} as summarised in the following steps.
	\begin{enumerate}
		\item[1.] Find points of ultimate tangentiality via~\eqref{thm1_ult_tan}.
		\item[2.] Solve for $\bar{d}$ via condition \eqref{eq_Hamil}. The solution may be a function of time, the state and/or adjoint.
		\item[3.] Integrate the system \eqref{sys_eq_1} and adjoint equation \eqref{eq_adj} backwards with $\bar{d}$ from the identified points $z$ to get candidate invariance barrier trajectories.
		\item[4.] Because the conditions of the theorem are necessary parts of come curves may have to be ignored.
\end{enumerate}
\section{Examples}\label{sec_examples}
\subsection{Constrained double integrator}

\begin{figure}[h]
	\begin{center}\
		\includegraphics[width=0.8\columnwidth]{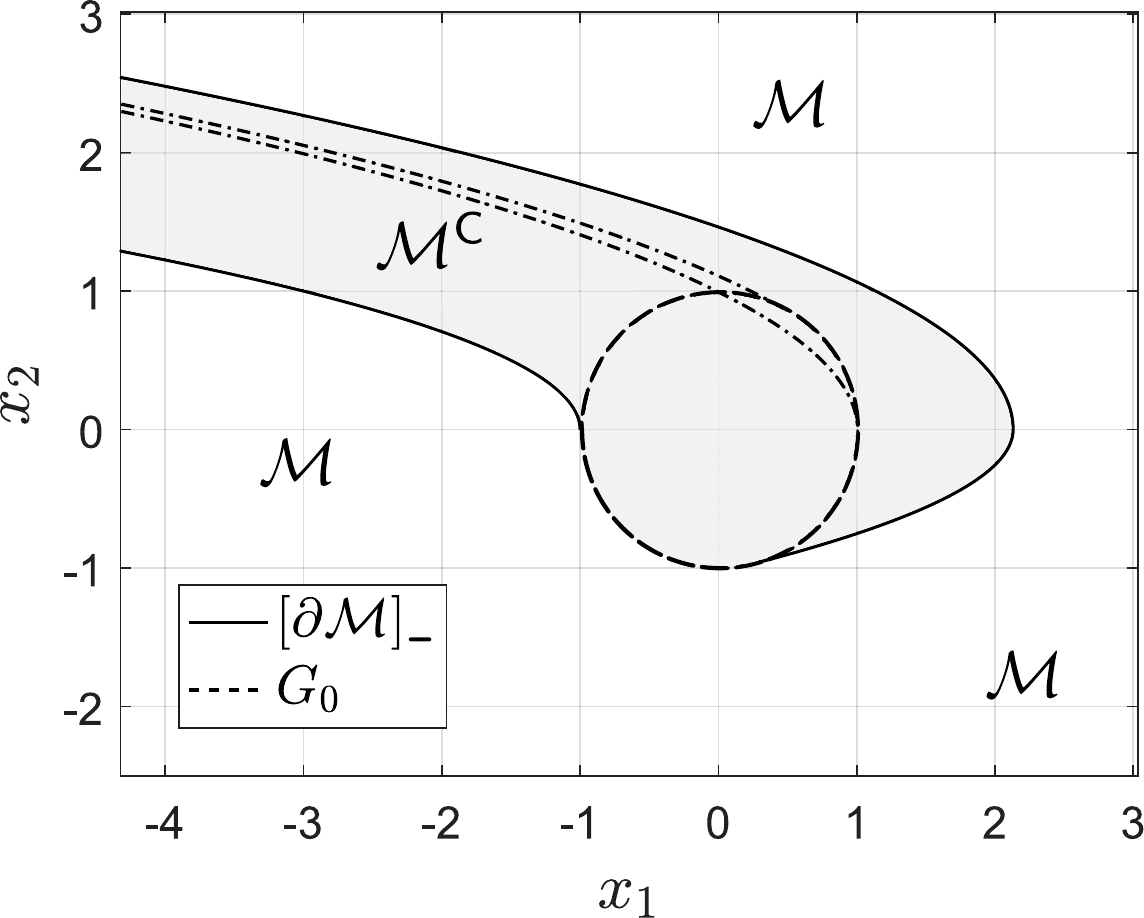} 
		\caption{The MRPI for the constrained double integrator.} 
		\label{fig:double_integrator2}                             
	\end{center}                                
\end{figure}

\begin{figure}[h]
	\begin{center}\
		\includegraphics[width=0.8\columnwidth]{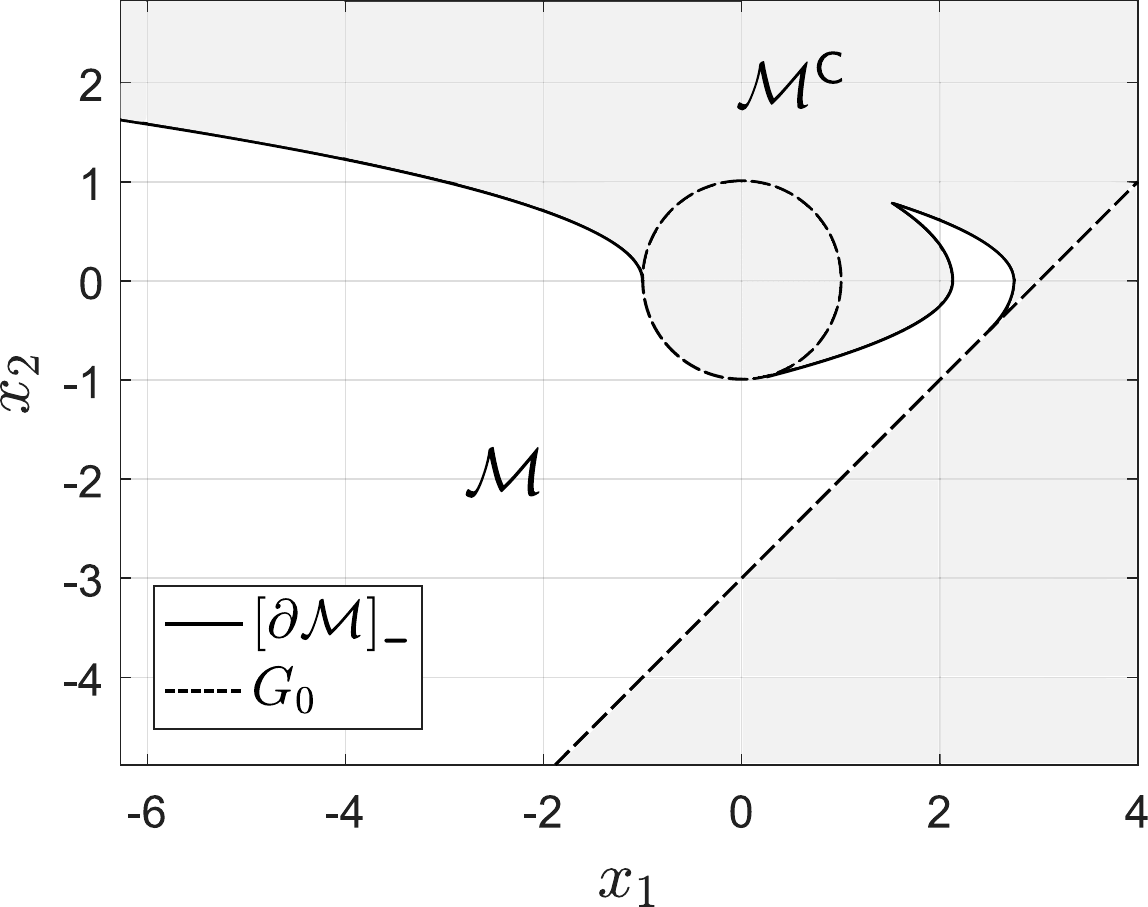} 
		\caption{The MRPI for the constrained double integrator with an additional linear constraint.} 
		\label{fig:double_integrator3}                             
	\end{center}                                
\end{figure}

\tolerance=2000
Consider the double integrator $\dot{x}_1 = x_2$, $\dot{x}_2 = d$, with $d\in[-0.5,-0.25]$ and $g_1(x) \triangleq -x_1^2 - x_2^2 + 1$. We invoke \eqref{thm1_ult_tan} to get:
	$
	\max_{d\in[-0.25,-0.5]} (-2x_1, -2x_2) (x_2, d)^T = 0,
	$
	from where we identify four points of ultimate tangentiality. The adjoint satisfies $\dot{\lambda}_1 = 0$, $\dot{\lambda_2} = -\lambda_1$ with $\lambda(\bar{t}) = (-2x_1(\bar{t}), -2x_2(\bar{t}))$, and from \eqref{eq_Hamil} we identify: $\bar{d}(t) = -0.25$ for $\lambda_2(t) \geq 0$; $\bar{d}(t) = -0.5$ for $\lambda_2(t) <0$. Integrating backwards, we find four candidate trajectories, as shown in Figure~\ref{fig:double_integrator2}. We ignore those drawn dash-dotted for otherwise $\DMO$ would include points on $G_0$ for which $\max_{d\in D} L_f g_1(z,d) > 0$.
We now add another constraint, $g_2(x) \triangleq x_1 - x_2 - 3$. We identify one point of ultimate tangentiality and find the candidate trajectory with the same $\bar{d}$. This curve intersects one of the curves associated with $g_1$ at a \emph{stopping point}, \cite{EL_IFACE2014}. We ignore the parts of both curves that extend, backwards in time, beyond this intersection point, because they are contained in parts of the state space for which either $g_1$ or $g_2$ may be violated by an admissible disturbance realisation. The set is shown in Figure~\ref{fig:double_integrator3}.

\subsection{Pendulum: nonlinear versus linearised model}
\begin{figure}[h]
	\begin{center}
		\includegraphics[width=0.8\columnwidth]{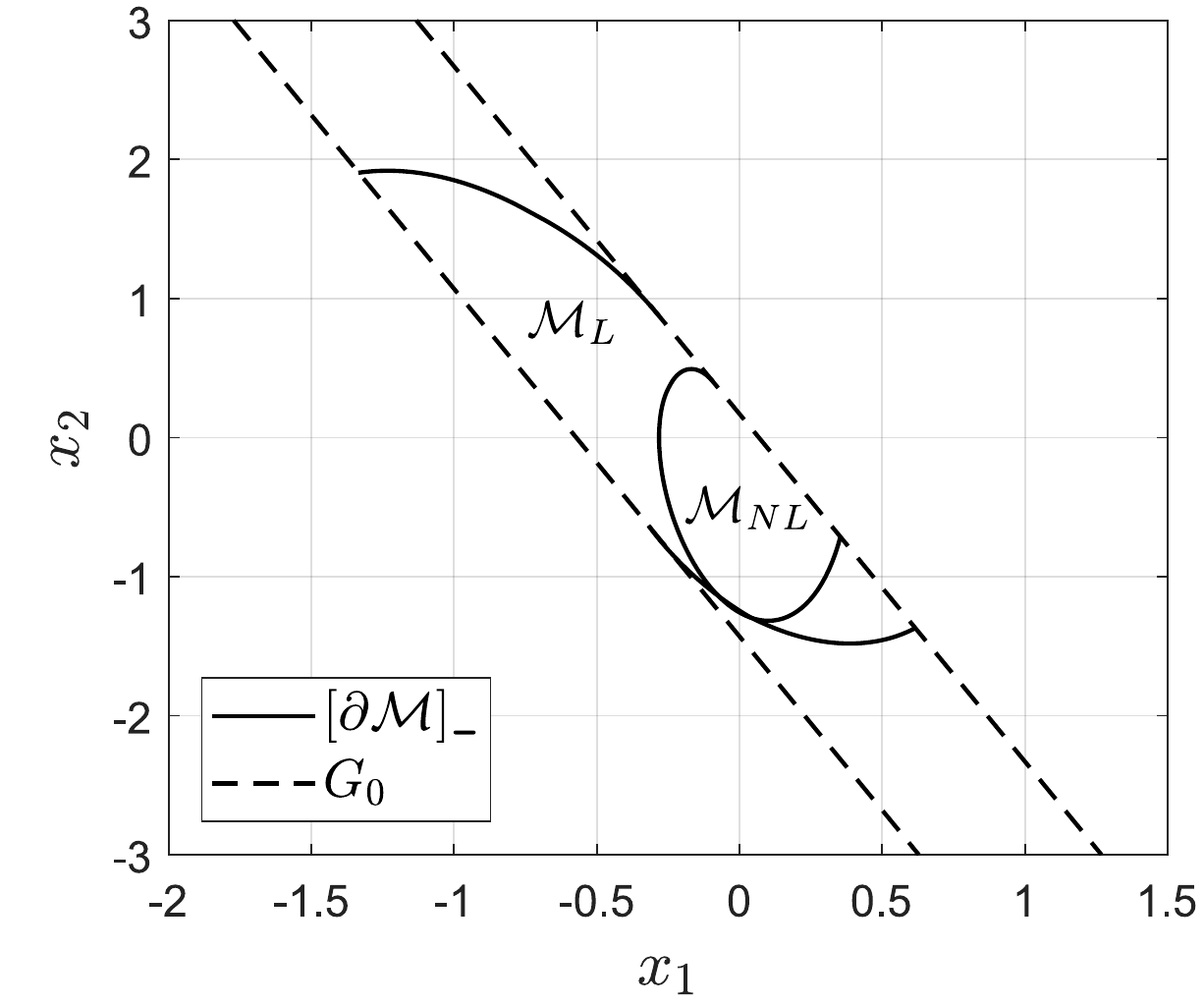}
		\vspace{-2mm}   
		\caption{The MRPIs for the nonlinear and linearised pendulum model ($\M_{NL}$ and $\M_{L}$ respectively), with $w = -0.3$.} 
		\label{fig:pendulum_lin}                             
	\end{center}                                
\end{figure}

We consider the nonlinear system of a pendulum actuated by a torque: $\ddot{\theta}(t) =-\frac{g}{l}\sin(\theta(t))+\frac{1}{ml^2} \tau(t) +d(t)$, with $\theta(t)$ the pendulum's angle; $\tau(t)$ the applied torque; $d(t)$ a disturbance torque; $g$ the gravitational constant; $l$ the pendulum's length; and $m$ the mass. Because $\sin(\theta) \leq 1$ for all $\theta$, (A3) is satisfied. It is desired that the actuator does not saturate during operation, hence, we impose $|\tau(t)|\leq2$. 
Assume that the disturbance is in the bounded interval $[-0.1,0.1]$. In \cite{Kolmanovsky_1998} the authors considered a linearised model and designed a linear feedback law, $\tau(\theta,\dot{\theta},w)=-k_1\theta-k_2\dot{\theta}+(k_1-1)w$, with the design parameters $k_1$ and $k_2$, and $w$ as constant that determines the equilibrium of $\theta$. We use $g=9.81$, $l=1$, $m=1$, $k_1=6.25$, $k_2=2.5$ and $w = -0.3$ and define $g_1(\theta, w)\triangleq\tau-2$ and $g_2(\theta) \triangleq -\tau - 2$. 
Invoking \eqref{thm1_ult_tan} we determine: $z\approx(-0.1, -1.17)$ and $z\approx(-0.09, 0.4)$. Integrating backwards we obtain $\M_{NL}$ as in Figure~\ref{fig:pendulum_lin}, where we have ignored the candidate invariance barrier trajectory ending at $(-0.1, -1.17)$. We also find the MRPI, labelled $\M_{L}$, for the linearised model considered in \cite{Kolmanovsky_1998}: $\dot{\theta}_1(t) = \theta_2(t)$, $\dot{\theta}_2(t) = \theta_1(t) + \tau(t) +d(t)$, with the same linear control law. We can see that the MRPI is in fact much smaller with the true nonlinear dynamics.
\vspace{-0.5cm}
\section{Conclusion}\label{sec_conclusion}
Extending results from \cite{DeDona_siam}, we have shown that parts of the boundary of the MRPI of a constrained nonlinear system are made up of integral curves that satisfy the maximum principle, and intersect the boundary of the constrained state space tangentiality. We used these facts to construct the set for some examples, illuminating some interesting properties.
\appendix
\vspace{-0.3cm}
\section{The Maximum Principle}
We reproduce the maximum principle, stated in terms of reachable sets, from \cite{Lee_Markus}.
\begin{theorem}[Maximum Principle]\label{thm_PMP}
	Consider \eqref{sys_eq_1} and $\bar{d}\in\D$ such that $x^{(\bar{d},x_0,t_0)}(t_1)\in\partial X_{t_1}(x_0)$ for some $t_1 > t_0$. Then, there exists a non-zero absolutely continuous maximal solution $\lambda^{\bar{d}}$ to the adjoint equation:
	\begin{align}
	&\dot{\lambda}^{\bar{d}}(t) = -\left( \frac{\partial f}{\partial x}(x^{(\bar{d},x_0,t_0)}(t),\bar{d}(t)) \right)^T \lambda^{\bar{d}}(t),\label{eq_1_PMP}
	\\
	&\mathrm{s. t.}\quad \max_{d\in D}\{\lambda^{\bar{d}}(t)^Tf(x^{(\bar{d},x_0,t_0)}(t),d) \}\nonumber \\
	& \qquad \;\;	= \lambda^{\bar{d}}(t)^T f(x^{(\bar{d},x_0,t_0)}(t),\bar{d}(t)) = \text{constant}\label{eq_2_PMP}
	\end{align}
	for almost every $t\in[t_0,t_1]$.
\end{theorem}

\bibliographystyle{IEEEtran}

\bibliography{bibliography_short}

\end{document}